\documentclass[8pt,reqno]{amsart}
\usepackage{latexsym}
\usepackage{amssymb}
\usepackage{amsmath}
\usepackage{amsthm}
\usepackage{nccmath}
\usepackage{amsfonts}
\usepackage{amssymb} 
\usepackage{booktabs}
\usepackage{rotating}
\usepackage{lscape}
\usepackage{color}

\usepackage{graphicx,color}  
\usepackage{graphpap, epic, xypic}  

\numberwithin{equation}{section}

\newcommand{\zz}{\pmb{\mathbb{Z}}}

\newcommand{\be}{\begin{equation}}
\newcommand{\ee}{\end{equation}}
\newcommand{\bes}{\begin{equation*}}
\newcommand{\ees}{\end{equation*}}

\newcommand{\eqn}{\begin{eqnarray}}
\newcommand{\feqn}{\end{eqnarray}}
\newcommand{\eqnn}{\begin{eqnarray*}}
\newcommand{\feqnn}{\end{eqnarray*}}

\newtheorem{theorem}{Theorem}

\newtheorem*{conj}{Conjecture-2.1"}

\DeclareMathOperator{\eup}{e}
\newcommand{\iup}{\textup{i}\hspace{1pt}}
\newcommand{\abs}[1]{\lvert#1\rvert}
\newcommand{\dup}{\textup{d}}

\begin{document}

\title{Compact Lie groups: Euler constructions and generalized Dyson conjecture.}
\author{Sergio L. Cacciatori}
\address{
Dipartimento di Scienza e Alta Tecnologia, Universit\`a degli Studi dell'Insubria,
Via Valleggio 11, 22100 Como, Italy, and INFN, via Celoria 16, 20133 Milano, Italy}
\email{sergio.cacciatori@uninsubria.it}
\author{Francesco Dalla Piazza}
\address{Dipartimento di Matematica, Universit\`a ``La Sapienza'', Piazzale A. Moro 2, I-00185, Roma, Italy}
\email{dallapiazza@mat.uniroma1.it, f.dallapiazza@gmail.com}
\author{Antonio Scotti}
\address{Dipartimento di Matematica,
Universit\`a degli Studi di Milano,  Via Saldini 50, 20133 Milano, Italy}
\email{ascotti@mindspring.com}


\begin{abstract}
A generalized Euler parameterization of a compact Lie group is a way for
parameterizing  the group starting from a maximal Lie subgroup, which allows a simple characterization of the range of parameters.
In the present paper we consider  the class of all compact connected Lie groups.
We present a general method for realizing their generalized Euler parameterization starting from any symmetrically embedded Lie group. Our construction is based on a detailed analysis of the geometry of these groups. As a byproduct  this gives rise to an interesting connection with certain Dyson integrals. In particular, we obtain a geometry based proof of a Macdonald conjecture regarding the Dyson integrals correspondent to the root systems associated to all irreducible symmetric spaces. As an application of our general method we explicitly parameterize all groups of the class of simple, simply connected compact Lie groups. We provide a table giving all necessary ingredients for all such Euler parameterizations.
\keywords{Lie groups \and Euler parameterization \and Macdonald conjecture \and Dyson integral}

\end{abstract}

\maketitle

\section{Introduction}\noindent
A simple procedure to 
parameterize compact simple Lie groups is provided by the generalized Euler parameterization
\cite{Cacciatori:2010dm,Bernardoni:2007rd,Bernardoni:2007rf,Cacciatori:2005yb,Cacciatori:2009qu,Cacciatori:2005gi,Bertini:2005rc,TS1,TS2}. 
The role of Lie groups both in mathematics and in physics is very important and well known. In several situations their local properties are sufficient, but there are also numerous concrete applications that require explicit realizations of the group matrices including the correct range of the parameters. This is the case in lattice gauge theories or in separability criteria for entangled configurations in quantum mechanics, just to cite two relevant physical examples.
For $SU(2)$ this problem has a natural solution in the Euler parameterization. A first important progress has been achieved in \cite{TS1} (followed by \cite{TS2}) where the authors have been able to define generalized Euler angles for $SU(N)$ groups. However, their construction was quite involved and does not allow an obvious generalization to other Lie groups. In \cite{Cacciatori:2005yb,Cacciatori:2005gi}, in studying the case of $G_2$, a more general strategy (later reviewed in \cite{Cacciatori:2009qu}), which was expected to be applicable case by case to all compact simply connected simple Lie groups, has been developed. Indeed, it has been first applied to the $SU(N)$ case in \cite{Bertini:2005rc}, showing the drastic simplification with respect to \cite{TS1,TS2}, and successively to the exceptional Lie groups $F_4$, $E_6$ and $E_7$ in \cite{Bernardoni:2007rf,Bernardoni:2007rd,Cacciatori:2010dm} respectively. It is only after these several years experience that we realized that the success of our method was due to the deep geometrical structure of the compact connected simple Lie groups, and, in fact, of all compact connected Lie groups. In the present paper we will analyze the geometry underlying our generalized Euler parameterization. This, beyond providing us with a final strategy valid for parameterizing all compact connected Lie groups, will also provide us with a direct connection between the geometry of compact Lie groups and certain integrals known as generalized Dyson integrals. As a consequence, we will obtain a geometrical proof of a class of particular cases of a conjecture stated in \cite{mac} by Macdonald, and already proved in general by Opdam in \cite{Opd}.

\subsection{Geometry of the Euler parameterization}
Let us consider the geometric structure underlying the generalized Euler parameterization of a simple Lie group.
Let $\pmb {\mathfrak {g}}\equiv {\rm Lie}(\pmb G)$ the Lie algebra of a real compact Lie group $\pmb G$. More precisely, we will assume that $\pmb {\mathfrak {g}}$ is
some matrix realization supporting a faithful representation of $\pmb G=\exp \pmb {\mathfrak{g}}$.
Our strategy is to start from a maximal symmetrically embedded proper subgroup $\pmb H$ of $\pmb G$. Let $\pmb G'$ be any real form of $\pmb G$ different from the compact one. Then $\pmb G'$ contains a maximal compact subgroup $\pmb H$ which is invariant under the action of the Cartan involution, and $\pmb H$ can be embedded in $\pmb G$. For any real form this is our choice for $\pmb H$ (in other words $\pmb G/\pmb H$ is the compact dual of the noncompact symmetric space $\pmb G'/\pmb H$). In particular when $\pmb G'$ is the split form we will call $\pmb H$ the
``the maximal compact subgroup of a split form $\pmb G'$ of $\pmb G$ (MCS)''.
This has the property ${\rm dim}\,\pmb G=2\dim\pmb H+r\equiv 2h+r$, where $r$ is the rank of $\pmb G$.
For the split case, the generalized Euler parameterization of $\pmb G$ takes the form
\begin{align}
\pmb G[x_1,\ldots, x_h; y_1,\ldots, y_r; z_1,\ldots, z_h]=&\pmb H[x_1,\ldots, x_h] \exp(y_1 c_1+\ldots+y_r c_r)\cdot \cr
&\cdot\pmb H[z_1,\ldots, z_h],
\end{align}
where $\pmb H[x_1,\ldots, x_m]$ is a parameterization of the maximal subgroup and $c_1,\ldots, c_r$ is a basis for a Cartan subalgebra $\pmb {\mathfrak{c}}$
in the complement of Lie$(\pmb H)$ in  $\pmb {{\mathfrak g}}$. Since the group is compact, one can choose $c_i$ so that the coordinates $y_i$ are periodic. A parameterization obtained in this way
in general is redundant for two reasons. \\
The first one is due to the fact that $\pmb H$ contains a finite subgroup $\pmb \Gamma$ of the maximal torus $\pmb T^r=\exp{\pmb {\mathfrak{c}}}$ of $\pmb G$.
Indeed, we will see that $\pmb \Gamma$ is isomorphic to $(\zz_2)^r$ if $\pmb G$ is simply connected, otherwise to one of its proper subgroups. Thus,
\begin{align*}
&\pmb H[x_1,\ldots, x_h] \exp(y_1 c_1+\ldots+y_r c_r) \pmb H[z_1,\ldots, z_h]=\pmb H[x_1,\ldots, x_h] \gamma^{-1}\cdot \cr
&\cdot \exp(y_1 c_1+\ldots+y_r c_r) \gamma
\pmb H[z_1,\ldots, z_h],
\end{align*}
for any $\gamma\in \pmb \Gamma$, so that we must reduce the range of the $\vec x$ coordinates w.r.t. the action of $\pmb \Gamma$:
$$
\pmb H[x_1,\ldots, x_h] \gamma^{-1}=\pmb H[\tilde x_1,\ldots, \tilde x_h].
$$
This is easily accomplished by accordingly restricting the range. \\
The second problem is due to the fact that the Weyl group $\pmb W$ acts nontrivially on $t\in \pmb T^r$,
$$
t\mapsto w^{-1}tw\in \pmb T^r,
$$
for any $w\in \pmb W$,
so also the range of the $\vec y$ coordinates has to be reduced.
We will show that this problem can be completely characterized in terms of the highest root
of $\pmb G$: after a suitable linear change of variables $s_i=\sum_{j=1}^r A_{ij} y^j$ we will
see that the right range of coordinates is expressed by the set of inequalities
\begin{eqnarray} \label{region}
0\leq n_1s_1+\ldots+n_r s_r \leq \pi, \qquad\ 0\leq s_i\leq \pi, \quad\ i=1,\ldots,r,
\end{eqnarray}
where $(n_1,\ldots,n_r)$ are the coefficients of the highest root $\tilde \alpha$ w.r.t. a basis of simple roots:
$\tilde \alpha=n_1 \alpha_1+\ldots+n_r \alpha_r$.
The volume of the whole group $\pmb G$ can be expressed in terms of the volume of the 
MCS subgroup $\pmb H$ times an integral directly connected to the generalized Dyson integral. On the other hand, the volumes of the compact Lie groups can be computed employing the Macdonald formula \cite{mac1}. Thus, incidentally, we see that our construction will turn out to be equivalent to prove certain particular cases of a conjecture due to Macdonald,
generalizing the Dyson integrals \cite{mac}.

This construction is more involved, but it works as well, for non simply connected compact Lie groups, as we will show in Section \ref{pigola}. Moreover, this parameterization applies to all compact connected Lie groups.
Furthermore, we will extend all the results to the case in which
a more general subgroup $\pmb H$ symmetrically
embedded in $\pmb G$ is considered, in place of the 
MCS one. In this case
the construction will turn out to be related to a version of Macdonald's conjecture for certain integrals associated
to non reduced root lattices.
In fact, the more interesting point is not
the proof of this conjecture, which can be obtained in a more general form using different methods (\cite{Opd}, see also \cite{OlPa}),
but its relation to the geometry of compact symmetric spaces.

We remark that this parameterization is also useful for concrete applications in Physics. Indeed, one often needs to work with an explicit realization of the parameterization of a Lie group,
including the right range for the parameters.


\subsection{Macdonald's conjecture}
We summarize the basic steps at the origin of  Macdonald's conjecture following the clear and punctual paper of P.~J.~Forrester and S.~O.~Warnaar \cite{ForWar}, to which we refer for a more extensive introduction.
The story of Macdonald's conjecture begins in the 1940s in the paper of Atle Selberg ``\"{U}ber einen Satz von A. Gelfond" \cite{sel1} where the author considered what is now known as Selberg integral:
\begin{align} \label{selbint}
S_n(\alpha,\beta,\gamma)&:=\int_0^1\cdots\int_0^1\prod_{i=1}^n t_i^{\alpha-1}(1-t_i)^{\beta-1}\prod_{1\leq i<j\leq n}|t_i-t_j|^{2\gamma}{\rm dt_1}\cdots{\rm dt_n} \cr
&=\prod_{j=0}^{n-1}\frac{\Gamma(\alpha+j\gamma)\Gamma(\beta+j\gamma)\Gamma(1+(j+1)\gamma)}{\Gamma(\alpha+\beta+(n+j-1)\gamma)\Gamma(1+\gamma)}.
\end{align}
This integral is valid for complex $\alpha$, $\beta$ and $\gamma$ such that:
\be
\Re(\alpha)>0, \qquad \Re(\beta)>0, \qquad \Re(\gamma)>-{\rm min}\{\frac 1n, \frac{\Re(\alpha)}{n-1},\frac{\Re(\beta)}{n-1}\},
\ee
corresponding to the domain of convergence of the integral. To contain the length of the paper, Selberg did not present there the proof of his claim, but in 1944, three years later, in the work ``Bemerkninger om et multiplet integral", \cite{sel2}. Note that the Euler beta integral is itself a Selberg integral with $n=1$.

For over thirty years, the Selberg integral was essentially unnoticed. The exception was a study by S.~Karlin and L.~S.~Shapley in 1953 \cite{karlshap} where they considered the special case $\alpha=1$, $\beta=1$ and $\gamma=2$ in relation to the volume of a certain moment space.
However, in the 1960s there were good reasons to make use of \eqref{selbint}. F.~J.~Dyson wrote a series of papers in the context of the statistical theory of energy levels of complex systems. A part of this series was written jointly with M.~L.~Mehta and published in 1963. Here, random Hermitian matrices were used to model highly excited states of complex nuclei. They considered systems with different symmetries described by matrices with real, complex or real quaternion elements. The ensembles of random matrices are called Gaussian orthogonal (GOE), unitary (GUE) and symplectic ensembles (GSE) respectively.
The joint probability density function for the three ensembles can be computed explicitly as:
\begin{equation}\label{fz}
\frac{1}{(2\pi)^{n/2} F_n(\beta/2)}
\prod_{i=1}^n {\rm e}^{-t_i^2/2}
\prod_{1 \le i < j \le n}
|t_i - t_j|^{\beta},
\end{equation}
where $\beta=1,2,4$ for the GOE, GUE and GSE respectively, and
$F_n$ is the normalization
\begin{equation}\label{fnz}
F_n(\gamma):=
\frac{1}{(2\pi)^{n/2}}
\int_{-\infty}^{\infty} \cdots \int_{-\infty}^{\infty}
\prod_{i=1}^n {\rm e}^{-t_i^2/2} \prod_{1 \le i < j \le n}
|t_i - t_j|^{2\gamma}\,{\rm d} t_1\cdots{\rm d} t_n,
\end{equation}
referred to as Mehta's integral, see \cite{mehta1}.
In \cite{MD} Mehta and Dyson evaluated $F_n(\beta/2)$ for each of
the three special values of $\beta$. Combining this with
the evaluations for $n=2$ and $n=3$ for general $\beta$ led them
to conjecture that
\begin{equation}\label{mehtacj}
F_n(\gamma) = \prod_{j=1}^n
\frac{\Gamma(1+j\gamma)}{\Gamma(1 + \gamma)}.
\end{equation}
The conjecture \eqref{mehtacj} can be proved evaluating the Mehta's integral using the Selberg integral, however in 1963 the Selberg's result was yet unknown.
The proof was finally given in the late 1970s by Enrico Bombieri.

The considerations on the symmetries of the complex systems that led to consider the three ensembles of Hermitian matrices can be applied also to unitary matrices
\cite{dyson1}. Making this choice of matrices, one obtains what are referred to as circular orthogonal ensemble (COE), circular unitary ensemble (CUE) and circular
symplectic ensemble (CSE). Their joint eigenvalues probability density function is given explicitly by:
\begin{equation}\label{Cp}
\frac{1}{(2\pi)^n C_n(\beta/2)} \prod_{1 \le i < j \le n}
\abs{\eup^{\iup \theta_i} - \eup^{\iup \theta_j} }^{\beta},
\end{equation}
where $C_n$ is the normalization
\begin{equation}\label{cnb}
C_n(\gamma) := \frac{1}{(2\pi)^n} \int_{-\pi}^{\pi} \cdots \int_{-\pi}^{\pi}
\prod_{1 \le i < j \le n}
\abs{\eup^{\iup \theta_i} - \eup^{\iup \theta_j} }^{2\gamma} \,
\dup \theta_1 \cdots \dup \theta_n,
\end{equation}
and $\beta = 1,2,4$ for the COE, CUE and CSE respectively.

As for \eqref{fnz}, the random matrix calculations give \eqref{cnb}
in terms of gamma functions for the three special values of $\beta$.
The case $n=2$ for general $\beta$ can be related to
the Euler beta integral, and the case $n=3$ gives
a sum which is a special instance of an identity of Dixon for a
well-poised $_3F_2$ series (cf. \cite{AAR99,ForWar} for details). Using these results, Dyson
made in  \cite{dyson1} the conjecture that:
\begin{equation}\label{CNb1}
C_n(\gamma)=\frac{\Gamma(1+n\gamma)}{\Gamma^n(1+\gamma)}.
\end{equation}
Moreover, Dyson observed that when $\gamma$ is a nonnegative
integer, say $k$,
\eqref{cnb} can be rewritten as the constant term (CT) in a Laurent
expansion. This allows \eqref{CNb1} to be rewritten as
\begin{equation}\label{CNb2}
\text{CT} \prod_{1 \le i < j \le n} \Bigl( 1 - \frac{x_i}{x_j} \Bigr)^k
 \Bigl( 1 - \frac{x_j}{x_i} \Bigr)^k = \frac{(kn)!}{(k!)^n}.
\end{equation}
This constant term identity and
the conjecture \eqref{CNb1}, were soon proved by J.~Gunson
and K.~Wilson \cite{Wilson62}, and later by
I.J.~Good \cite{Good70}.
R.~Askey \cite{Askey80} observed that
the Selberg integral can be used to prove Dyson's conjecture
\eqref{CNb1} directly.

The Macdonald's conjecture \cite{mac} may be considered as a generalization of the Dyson's conjecture \eqref{CNb2}. Let $R$ be a reduced root system, $e^\alpha$
denote the formal exponential corresponding to $\alpha\in R$ and $k$ a nonnegative integer. Then Macdonald conjectured (cf. \cite{mac}, Conjecture 2.1) that the
constant term in the polynomial
\be
\prod_{\alpha\in R}(1-e^\alpha)^k
\ee
should be equal to $\prod_{i=1}^l\binom{kd_i}{k}$, where the $d_i$ are the degrees of the fundamental invariants of the Weyl group of $R$ and $l$ the rank of $R$.
Macdonald wrote this relation in an equivalent form which will turn out to be useful later. Let $G$ be a compact connected Lie group, $T$ a maximal torus of $G$,
such that $R$ is the root system of $(G,T)$
and define:
\be
\Delta(t)=\prod_{\alpha\in R^+}(e^{\alpha/2}(t)-e^{-\alpha/2}(t)),
\ee
where
$t\in T$,
the exponentials are regarded as characters of $T$  and $R^+$ is a choice of positive roots. Then $|\Delta(t)|^2=\prod_{\alpha\in R}(1-e^\alpha(t))$ is a positive
real-valued continuous function on $T$. This function enters in Weyl's integration formula
\be \label{intfGT}
\int_G f(x)dx=\frac{1}{|W|}\int_T |\Delta(t)|^2f(t)dt
\ee
for any continuous class function $f$ on $G$. In \eqref{intfGT}, $dx$ and $dt$ are the normalized Haar measure on $G$ and $T$ respectively
$(\int_G dx=\int _T dt=1)$. Thus, the conjecture can be rewritten as (cf. \cite{mac} Conjecture 2.1'):
\be \label{cj2.1'}
\int_T|\Delta(t)|^{2k}dt=\prod_{i=1}^l\binom{kd_i}{k}.
\ee
The equivalence of the two formulations follows from the fact that the integration over $T$ kills all but the trivial character, or in other words, selects the
constant term in $|\Delta(t)|^{2k}=\prod_{\alpha\in R}(1-e^\alpha(t))^k$. An observation that generalizes further the conjecture is that \eqref{cj2.1'} makes
sense if the integer $k$ is replaced by a complex number, $s$, with positive real part, $\Re(s)>0$. In this case the right hand side is replaced by
\be
\prod_{i=1}^l\frac{\Gamma(sd_i+1)}{\Gamma(s+1)\Gamma(sd_i-s+1)}.
\ee

In the same paper Macdonald generalized the conjecture further (cf. \cite{mac}, Conjecture 2.3). For this, let $R$ be a root system, now not necessarily reduced,
and for each $\alpha\in R$ let $k_\alpha$ be a nonnegative integer such that $k_\alpha=k_\beta$ if $|\alpha|=|\beta|$, then the constant term in the Laurent
polynomial
\be
\prod_{\alpha\in R}(1-e^\alpha)^{k_\alpha}
\ee
should be equal to the product
\be
\prod_{\alpha\in R}\frac{\left(\left|\langle\rho_k,\check \alpha \rangle+k_\alpha+\frac 12 k_{\alpha/2}\right|\right)!}{\left(\left|\langle\rho_k,\check \alpha
\rangle+\frac 12k_{\alpha/2}\right|\right)!},
\ee
where $\rho_k=\frac 12\sum_{\alpha\in R^+}k_\alpha \alpha$, $\check\alpha=\frac{2\alpha}{|\alpha|^2}$ is the coroot corresponding to $\alpha$, $k_{\alpha/2}=0$
if $\frac 12\alpha \not\in R$ and $\langle,\rangle$ is the usual scalar product induced by the Killing form. When the $k_\alpha$ are all equal this reduces to
the previous conjecture.

The Macdonald's conjecture was finally proved in a slightly more general form by Opdam \cite{Opd} considering $k_\alpha$ a complex valued Weyl invariant function
with positive real part. This is the content of Theorem 4.1 of \cite{Opd}:
\begin{theorem}[Macdonald-Opdam] \label{thM-O}
Let $R$ be a possibly non reduced root system, and let $k\in \mathbb{K}$ such that\footnote{The vector space $\mathbb{K}\simeq \mathbb{C}^n$ is the space of all
complex valued Weyl invariant functions on $R$, $m$ equals the numbers of conjugacy classes of roots in $R$ and elements of $\mathbb{K}$ are called multiplicity
functions on $R$. The notation $k_\alpha$ denote the evaluation of $k\in \mathbb{K}$ on $\alpha\in R$.} $\Re(k_\alpha)\geq 0$, $\forall \alpha\in R$. Then
\be
\int_T\sigma(k,t)dt=\prod_{\alpha\in R^+}\frac{\Gamma(\langle\rho(k),\check{\alpha}\rangle+k_\alpha+\frac 12 k_{\alpha/2}+1)\Gamma(\langle\rho(k),
\check{\alpha}\rangle-k_\alpha-\frac 12 k_{\alpha/2}+1)}{\Gamma(\langle\rho(k),\check{\alpha}\rangle+\frac 12 k_{\alpha/2}+1)\Gamma(\langle\rho(k),
\check{\alpha}\rangle-\frac 12 k_{\alpha/2}+1)},
\ee
\noindent
where $\sigma(k,t)=\prod_{\alpha\in R^+}|t^{\frac\alpha 2}-t^{-\frac\alpha 2}|^{2k_\alpha}$ and $T$ is the compact part in the ``polar decomposition" of the
maximal torus.
\end{theorem}


\section{The split case}\label{sec:mms}
In this section we describe our general construction of the Euler parameterization for an arbitrary compact connected Lie group, $\pmb G$ relative to a specific choice of the corresponding subgroups $\pmb H$.
\subsection{The compact Lie groups}\label{pippo}
\noindent
Let $\pmb G_0$ be a real compact connected semisimple Lie group (cf. the definition in \cite{Hel}, page 131). This means that
$\pmb G_0\simeq \pmb G_1 \times \pmb G_2 \times \dots \times \pmb G_n$, where $\pmb G_i$, $i=1,\ldots n$ are simple Lie groups uniquely determined
(up to permutations).
Let $\pmb H\simeq \pmb H_1 \times \pmb H_2 \times \dots \times \pmb H_n$ 
such that
$\pmb H_i$ is a maximal Lie subgroup symmetrically embedded in $\pmb G_i$. Then $\pmb H$ is connected (see \cite{Hel}, Chapter VI, Theorem 1.1). In general, $\pmb H_i$ is not simple, nor semisimple, but it has the form
$\pmb H_i \simeq \pmb H_{0,i}\times \pmb T^{s_i}/\Delta_i$, where $\pmb H_{0,i}$ is semisimple, $\pmb T^{s_i}$ is an Abelian torus, and $\Delta_i$ a finite
subgroup. Since our aim is to construct the Euler parameterization of $\pmb G_0$ relative to the subgroup $\pmb H$, and then applying the same procedure to $\pmb H$ inductively, we are forced to
consider
the more general case
\begin{eqnarray} \label{GG0TD}
  \pmb G \simeq \pmb G_0 \times \pmb T^s/ \Delta,
\end{eqnarray}
with $\pmb G_0$ as before, $\pmb T^s$ an Abelian torus, and $\Delta$ a finite subgroup. 
\vskip 0.25cm
\noindent
{\bf Remark.} 
The class of compact connected Lie groups of the form \eqref{GG0TD} coincides with the class of all compact connected Lie groups. We need only to show that any compact connected Lie group $\pmb G$ has the form \eqref{GG0TD}.  Let $\pmb G_0 := \pmb G'$ be the derived group and let the torus $T^s =\pmb Z^0$ be the connected component of the identity of the center $\pmb Z$ of $\pmb G$. Then, the multiplication map $m:\pmb G'\times \pmb Z^0 \rightarrow \pmb G$ is surjective. Indeed, it is a homomorphism since $\pmb Z^0$ is central. Moreover, since $\rm{Lie}\,(\pmb G^0)\oplus \rm{Lie}\, (\pmb Z^0)=\rm{Lie}\,(\pmb G)$, the differential of $m$ is surjective, so the image of $m$ is open. Since $\pmb G$ is connected, it follows that $m$ is surjective. The kernel of $m$ is obtained by embedding $\pmb G' \cap \pmb Z^0$ in $\pmb G'\times \pmb Z^0$ via $\gamma\rightarrow (\gamma,\gamma^{-1})$, $\gamma \in \pmb G'\cap \pmb Z^0$. 
The image of such map is the kernel of $m$ and is a finite group.
\vskip 0.25cm
The parameterization of $\pmb G$ from $\pmb G_0$ is quite elementary
and we can concentrate here on the parameterization of $\pmb G_0$ only. From now on we will assume
\be
\pmb G\equiv \pmb G_0.
\ee
Note that $\pmb H$ is symmetrically embedded
in $\pmb G$, but is not maximal, unless $\pmb G$ is simple. 

Let $\pmb{\mathfrak g}$ and $\pmb{\mathfrak h}$ be the Lie algebras associated to $\pmb G$ and $\pmb H$ respectively.
In this section, we will assume that $\pmb H$ is
MCS:
with this, we mean that $\pmb H_i$ is a 
MCS
subgroup of $\pmb G_i$.
Since in this case ${\rm rank}\, (\pmb{G/H})={\rm rank}\, \pmb G$, we can choose a Cartan subalgebra $\pmb{\mathfrak{c}}$ of $\pmb{\mathfrak g}$ such that $\pmb{ \mathfrak{c}}\cap \pmb{\mathfrak{h}}=0$.
Thus, the generalized Euler parameterization of the group $\pmb G$ w.r.t. $\pmb H$ takes the form
\begin{eqnarray}
\pmb{G}=(\pmb{H'}/\pmb \Gamma) e^{\pmb{\mathfrak{c}}} \pmb H,\label{param}
\end{eqnarray}
where $\pmb{H'}$ is a copy of $\pmb H$ and $\pmb \Gamma=\pmb H\cap e^{\pmb{\mathfrak{c}}}$ is a finite subgroup of the maximal torus that will be specified
later. In Section \ref{sec:gms} we will extend the parameterization to an arbitrary symmetrically embedded subgroup.
Before entering the details of the construction, we need to specify some further technical facts.
\subsection{Some technical facts and definitions} \label{pigola}
Any finite dimensional semisimple Lie algebra $\pmb{\mathfrak{g}}$ admits a unique compact form. There is a unique (up to isomorphisms) simply connected,
compact Lie group $\tilde {\pmb G}$ having $\pmb{\mathfrak{g}}$ as the associated Lie algebra. However, more in general, there are more than one connected compact Lie
group having the same Lie algebra. These are a finite number
\begin{eqnarray}
  \pmb G^k=\tilde {\pmb G}/\Gamma^k, \qquad\ k=1,2,\ldots,m,
\end{eqnarray}
where $\Gamma^k$ are finite subgroups of the center $Z$ of $\tilde{\pmb G}$. Notice that $\tilde {\pmb G}$ has finite center (cf. \cite{bump}, Proposition 23.11, pag. 200).
In particular, we
set
$\pmb G^1\equiv \tilde{\pmb G}$ and $\pmb G^m\equiv \tilde{\pmb G}/Z=: \pmb G_Z$. Then each $\pmb G^k$ is a covering of $\pmb G_Z$
and is covered by $\tilde {\pmb G}$. It is known that each of such groups admits a faithful linear representation (cf. \cite{bump}, Theorem 4.2, pag. 26). Let $(R_k,V_k)$ such a representation
for $\pmb G^k$ (in particular $(Ad,\pmb{\mathfrak {g}})$ is faithful for $\pmb G_Z$). It induces a faithful representation $(\rho_i,V_i)$ of
$\pmb{\mathfrak{g}}$, so that the following diagram is commutative:

\

\hspace{3cm}
{\centering
\xymatrix{
{\pmb G^k} \ar@{<-}[d]_{\exp_{\pmb G^k}}
&Aut(V_k) \ar@{<-}[l]_{\!\! {R_k} } \ar@{<-}[d]^{\rm Exp} \\
{\pmb {\mathfrak {g}}} \ar@{->}[r]^{\rho_k} & End(V_k)
}
}\\[.5cm]
Since $\pmb G^k$ is compact, $R_k$ is injective and continuous, and ${\rm Aut}(V_k)$ is $T_2$, then $\pmb G^k$ and $R_{k}(\pmb G^k)$ are homeomorphic
(see \cite{kos}, Theorem 8.8) and, in particular, they have the same fundamental group\footnote{We are grateful to S. Pigola for explaining us these points.}. More in general, this means that we can construct a realization of
the desired compact form simply by exponentiating the matrices associated to the Lie algebra $\pmb {\mathfrak {g}}$ via the representation $\rho_k$
induced by the faithful representation $R_k$ of $\pmb G^k$. For this reason, we will call $\rho_k$ a {\it $\pmb G^k$-faithful representation}. Thus,
we parameterize the desired compact $\pmb G$ form by working with the right $\pmb G$-faithful representation of the algebra.\\

Let $\pmb H$ be a subgroup of $\pmb G$ as defined previously. If we are working with a $\pmb G$-faithful representation $(\rho, V)$ for $\pmb {\mathfrak g}$,
then $\rho$ will decompose into a direct sum of representations of the Lie algebra $\pmb {\mathfrak h}$ of $\pmb H$, among which at least one is surely
$\pmb H$-faithful (whereas the complementary ones will give rise to compact forms of the group covered by $\pmb H$). As a consequence, one can construct
the corresponding parameterization of $\pmb H$ by worrying about the $\pmb H$-faithful representation only, which is automatically present in the decomposition.

\subsection{Parameterization}
The problem of parameterizing $\pmb H$ and $\pmb {H'}$ is then the same as for $\pmb G$ and can be obtained inductively. Thus, if we want to get an almost
everywhere one to one parameterization of $\pmb{G}$, the only problem is to determine the right range for the parameterization of the toric part
\begin{eqnarray}
e^{\pmb{\mathfrak{c}}}[y_1,\ldots, y_r]=\exp(y_1 c_1+\ldots+y_r c_r).
\end{eqnarray}
Now, we will see that the range for the $y$'s is independent from the starting $\pmb G$-faithful representation, depending on the Adjoint representation only.
This means that such determination is in a sense universal, and the details discriminating among the different compact forms of $\pmb G$ will depend only
from the periodicities of the $U(1)$ factors entering the parameterizations, and the action of the finite subgroups. Furthermore, employing the isomorphism $\pmb G \simeq \pmb G_1 \times \pmb G_2 \times \dots \times \pmb G_n$, we can parameterize each factor independently. Thus without loss of generality, we focus on the case when $\pmb G$ is simple. In particular, then, $\pmb H$ is maximal in $\pmb G$.\\
Using the notation in Section \ref{pippo}, we write $\pmb{\mathfrak {g}}=\pmb{\mathfrak {h}}\oplus \pmb{\mathfrak {p}}=\pmb{\mathfrak {h}}\oplus \pmb{\mathfrak{c}}\oplus \pmb{\mathfrak {p}'}$, with
$\dim\pmb{\mathfrak {h}}=\dim\pmb{\mathfrak {p}'}= h$ and $\dim \pmb{\mathfrak{c}}=r$. Now consider the complexification $\pmb{\mathfrak {g}}_{\pmb{\mathbb C}}$ of $\pmb{\mathfrak{g}}$. It has also the
decomposition $\pmb{\mathfrak {g}}_{\pmb{\mathbb C}}=\pmb W_- \oplus \pmb{\mathfrak{c}} \oplus \pmb W_+$, where $\pmb W_\pm$ is the direct sum of the root spaces
$\pmb W_{\pm\alpha}$ such that $\alpha$
is a positive root. We can thus pick out the following two bases for $\pmb{\mathfrak {g}}_{\pmb{\mathbb C}}$:
\begin{itemize}
\item $\{c_i\}_{i=1}^r \cup \{\lambda_{\alpha_a} \cup \lambda_{-\alpha_a}\}_{a=1}^h,$ where $\lambda_{\alpha}$ is the eigenvector corresponding
to the root $\alpha$ and the $\alpha_a$ are the positive roots;
\item  $\{c_i\}_{i=1}^r \cup \{ t_a \}_{a=1}^h \cup \{ p_b\}_{b=1}^h$, where $t_a$ and $p_b$ generate $\pmb{\mathfrak h}$ and
$\pmb{\mathfrak {p}'}$
respectively, and are chosen so that ${\rm ad}_{t_a}$ and ${\rm ad}_{p_b}$ are diagonalizable, and the decomposition is Killing orthogonal.
\end{itemize}
Notice that only the second one is a basis for the compact algebra $\pmb{\mathfrak g}$. It satisfies the following relations:
\begin{eqnarray}
&& [t_a, t_b]\in \pmb{\mathfrak {h}}, \qquad [t_a, c_i]\in \pmb{\mathfrak {p}'}, \qquad [t_a, p_b]\in \pmb{\mathfrak {p}},\cr
&& [c_i, c_j]=0, \qquad [c_i, p_b]\in \pmb{\mathfrak {h}}, \qquad [p_a, p_b]\in \pmb{\mathfrak {h}}.
\end{eqnarray}
Indeed, the maximal symmetrically embedded compact subalgebras are in biunivocal correspondence with the real forms or, equivalently, with the Cartan decompositions of the algebra.
This means that there exists an involution $\theta: \pmb{\mathfrak {g}} \rightarrow \pmb{\mathfrak {g}}$,
$\theta^2={\rm id}$, such that $\pmb{\mathfrak h}$ and $\pmb{\mathfrak p}$
are the corresponding eigenspaces, with eigenvalues $1$ and $-1$ respectively. Since $\theta$ is a homomorphism, this implies
\begin{eqnarray}
[\pmb{\mathfrak {h}}, \pmb{\mathfrak {h}}]\subseteq \pmb{\mathfrak {h}},
\qquad [\pmb{\mathfrak {h}}, \pmb{\mathfrak {p}}]\subseteq \pmb{\mathfrak {p}},
\qquad [\pmb{\mathfrak {p}}, \pmb{\mathfrak {p}}]\subseteq \pmb{\mathfrak {h}},
\end{eqnarray}
and the Killing orthogonality between the two spaces. Moreover, ${\rm ad}$-invariance of the Killing form $\langle, \rangle$ implies
\begin{eqnarray}
\langle [t_a, c_i], c_j  \rangle=\langle t_a, [c_i, c_j]  \rangle=0,
\end{eqnarray}
so that $[t_a, c_i]\in \pmb{\mathfrak {p}'}$.\\
The rules $[c_i, p_b]\in \pmb{\mathfrak {h}}$, $[t_a,c_i]\in \pmb{\mathfrak {p'}}$ allow to provide a simple relation between
the two bases defined above. Indeed, we can obtain from $\pmb{\mathfrak g}$ a new real form of $\pmb{\mathfrak g}_{\pmb{\mathbb C}}$ by means of the Weyl's unitary
trick which consists in defining the new generators
\begin{itemize}
\item $\tilde t_a= t_a, \qquad \tilde c_j=i c_j, \qquad \tilde p_b=ip_b.$
\end{itemize}
This is the noncompact form $\pmb{\mathfrak {g}}_{(r)}$, with signature $r$. In this case, the operators ${\rm ad}_{\tilde h_j}$ are represented by
symmetric matrices since ${\rm ad}$-invariance and symmetry of the Killing form give
\begin{eqnarray}
\langle [\tilde c_i,\tilde p_a], \tilde t_b \rangle =-\langle \tilde p_a, [\tilde c_i, \tilde t_b] \rangle,
\end{eqnarray}
and the form is positive definite over the $\tilde t_b$ and negative over the complementary space. This means that such matrices can be
diagonalized, with real eigenvalues, by means of real combinations of the vectors $\tilde t_a$, $\tilde p_b$. Then, an eigenvector
corresponding to a non-zero root $\alpha$ will have the form $\lambda_\alpha= t_\alpha +i p_\alpha$, with $t_\alpha \in \pmb{\mathfrak h}$ and
$p_\alpha\in \pmb{\mathfrak p}$. Notice that both $t_\alpha$ and $p_\alpha$ are necessarily non vanishing. Indeed,
$[c,\lambda_\alpha]=\alpha(c) \lambda_\alpha$ for all $c\in \pmb{\mathfrak{c}}$ implies
\begin{eqnarray}
[c, t_\alpha]=i \alpha(c) p_\alpha, \quad [c, p_\alpha]=-i \alpha(c) t_\alpha,
\end{eqnarray}
and $t_\alpha=0$ or $p_\alpha=0$ would imply $\alpha(c)=0$ for any $c\in \pmb{\mathfrak{c}}$.\\
In conclusion, we can choose the basis
$t_a$, $p_b$, $c_i$, so that the relation between the two bases is
\begin{eqnarray}
\lambda_{\alpha_a}= t_a +i p_a, \qquad \lambda_{-\alpha_a}= t_a -i p_a, \qquad a=1,\ldots,h.
\end{eqnarray}
Moreover, since $\langle \lambda_\alpha, \lambda_\beta \rangle\neq 0$ if and only if $\alpha+\beta=0$, we can normalize the basis so that
it becomes an orthonormal basis.\\
Now, we show that these known facts have interesting consequences for the Euler parameterization. We can write
\be
\pmb G[x_1,\ldots, x_h, y_1, \ldots, y_r, z_1, \ldots, z_h]= e^{\sum_{a=1}^h x_a t_a} e^{\sum_{i=1}^r y_i c_i} e^{\sum_{b=1}^h z_b t_b}
\equiv
(\pmb{H'}/\pmb \Gamma ) e^{\pmb{\mathfrak{c}}} \pmb H.
\ee
The invariant measure expressed in terms of this parameterization is
\begin{eqnarray}
d\mu_{\pmb G}[\vec x; \vec y; \vec z]=d\mu_{\pmb H} [\vec z] d\mu_{\pmb B} [\vec x; \vec y],
\end{eqnarray}
where $d\mu_{\pmb H} [\vec z]$ is the invariant measure associated to $\pmb H$ and
\begin{eqnarray}
d\mu_{\pmb B} [\vec x; \vec y]= \det {J} (\vec x, \vec y) \prod_{a=1}^h dx_a \prod_{i=1}^r  dy_i,
\end{eqnarray}
$J$ being the $h\times h$ matrix with components
\begin{eqnarray}
J^a_{\ b}:= \langle e^{-\pmb{\mathfrak{c}}} {\pmb{H'}}^{-1} \frac {\partial \pmb{H'}}{\partial x_a} e^{\pmb{\mathfrak{c}}}, p_b \rangle.
\end{eqnarray}
Notice that ${\pmb{H'}}^{-1} d\pmb{H'} =:J_{\pmb H}$ is the left invariant one form for the $\pmb H$ subgroup in the $\pmb{H'}$ parameterization,
$J_{\pmb H}= \sum_{a=1}^h J_{\pmb H}^a t_a$. Thus
\begin{eqnarray}
d\mu_{\pmb B} [\vec x; \vec y]= d\mu_{\pmb H} [\vec x] \det M \prod_{i=1}^r  dy_i, \qquad M^a_{\ b}:= \langle e^{-\pmb{\mathfrak{c}}} t_a e^{\pmb{\mathfrak{c}}}, p_b \rangle.
\end{eqnarray}
Now $t_a= (\lambda_{\alpha_a}+\lambda_{-\alpha_a})/2$ so that
\begin{eqnarray} \label{eq:adj}
e^{-\pmb{\mathfrak{c}}} t_a e^{\pmb{\mathfrak{c}}}=
\cosh (\alpha_a (\pmb{\mathfrak{c}})) t_a +i \sinh (\alpha_a (\pmb{\mathfrak{c}})) p_a.
\end{eqnarray}
Since the roots are real on $\tilde c_i$, if we define $\vec \alpha_a\equiv (\alpha_a^1,\ldots, \alpha_a^r)$
with $\alpha_a^i=\alpha_a(\tilde c_i)$, we get $\alpha_a (\pmb{\mathfrak{c}})=-i\sum_{i=1}^r \alpha_a^i y_i\equiv -i\vec \alpha_a \cdot \vec y$. Then
\begin{eqnarray}
\det M=\prod_{a=1}^h \sin (\vec \alpha_a \cdot \vec y).
\end{eqnarray}
Thus, the invariant measure takes the form
\begin{eqnarray} \label{eq:misure}
d\mu_{\pmb G}[\vec x; \vec y; \vec z]=d\mu_{\pmb H} [\vec z] d\mu_{\pmb H} [\vec x] \prod_{a=1}^h \sin (\vec \alpha_a \cdot \vec y)
\prod_{i=1}^r  dy_i.
\end{eqnarray}
The range of the $z$ coordinates is such to cover the subgroup $\pmb H$, whereas the range $R_y$ for the $y$ coordinates is defined by
the conditions $0\leq \vec \alpha_a \cdot \vec y \leq \pi$, and the range for the $x$ coordinates is such to cover $\pmb{H'}/{\pmb \Gamma}$. In particular, as a consequence of equation \eqref{eq:adj}, the range for the $y_i$'s depends on the adjoint representation and not on the particular $\pmb G$-faithful representation we are considering.
Notice that equation \eqref{eq:misure} implies the interesting relation
\begin{eqnarray}\label{integral}
\int_{R_y} \prod_{a=1}^h \sin (\vec \alpha_a \cdot \vec y) \prod_{i=1}^r  dy_i=\frac {{\rm Vol} (\pmb G)\ |\pmb\Gamma|}{{\rm Vol} (\pmb H)^2},
\end{eqnarray}
where the volumes can be computed by means of the Macdonald's formula \cite{mac1} and $|\pmb\Gamma|$ is the cardinality of $\pmb\Gamma$.\\
When $\pmb G$ is simply connected, it is easy to see that ${\pmb\Gamma} \simeq \zz_2^r$. Indeed, the elements of ${\pmb \Gamma}={\pmb H}\cap e^\mathfrak{{\pmb c}}$
are the elements of $e^\mathfrak{{\pmb c}}$ whose square is the identity (see \cite{Hel}, section VII, Theorem 8.5). Since the basis $c_1,\ldots, c_r$ of
$\mathfrak{{\pmb c}}$ can be chosen so that $e^{tc_i}$ has period $T$, $\pmb\Gamma$ is generated by $e^{\frac T2 c_i}$ that proves our claim.
In particular, $|\pmb\Gamma|=2^r$.\\
When $\pmb G$ is not simply connected, this is not true in general and $\Gamma$ is isomorphic to a proper subgroup of $\zz_2^r$. Indeed, if
$\Phi: \tilde {\pmb G} \longrightarrow \pmb G$ is the universal covering map, then $\Gamma \simeq \Phi(\zz_2^r)$.

\subsection{Connections with the generalized Dyson integrals}\label{sec:dyson}
Let us now look closer at the integral (\ref{integral}). It is convenient to introduce the following change of variables. Let
$\alpha_{a_1}, \ldots, \alpha_{a_r}$ be simple roots. Then, we define the new coordinates $s_i$, $i=1,\ldots,r$ by
\begin{eqnarray}
s_i:=\vec y \cdot \vec \alpha_{a_i}.
\end{eqnarray}
From this we get
\begin{eqnarray}
ds_1 \wedge \ldots \wedge ds_r= V_F \prod_{i=1}^r \frac {\|\alpha_{a_i}\|^2}2 dy_1 \wedge \ldots \wedge dy_r ,
\end{eqnarray}
where $V_F$ is the volume of the fundamental region (parallelogram) defined by the simple coroots $\check \alpha_{a_i}$.
Then, the integral in \eqref{integral} takes the form
\begin{eqnarray}
I=\frac {2^r}{V_F \prod_{i=1}^r \|\alpha_{a_i}\|^2} \int_{\tilde R_s}
\prod_{a=1}^h \sin (\vec n_a \cdot \vec s) \prod_{i=1}^r  ds_i,
\end{eqnarray}
where $\vec n_a$ are the coordinates of the positive roots expressed w.r.t. the simple roots and take value in $\pmb{\mathbb {N}}^r$ and $\tilde R_s$ is the range for the $s$ coordinates.
In particular, for the simple roots we have
\begin{eqnarray}
\prod_{i=1}^r \sin (\vec n_{a_i} \cdot \vec s)=\prod_{i=1}^r \sin (s_i),
\end{eqnarray}
so that the range of coordinates is a subset of the cube $0\leq s_i\leq \pi$.
The remaining conditions are $0\leq \vec n_{a} \cdot \vec s\leq \pi$ for all the other positive roots.
The hyperplanes $\vec n_{a} \cdot \vec s=k\pi$, with $k$ integer, cut the cube in a tiling whose sectors are all equivalent being related each others by the Weyl reflections.
We know that for a simple group the highest root (relative to the simple root system) $\alpha_{\tilde a}=\sum_{i=1}^r \tilde n_i \alpha_{a_i}$
has the property
$\tilde n_i \geq n_a^i$ for all $a$ and $i$ (indeed, it is nothing but the highest weight of the adjoint representation).
Thus, the inequalities $0\leq \vec n_{a} \cdot \vec s\leq \pi$ defining the tiling reduce just to one. Indeed,
\begin{eqnarray}
0\leq \vec {\tilde n} \cdot \vec s\leq \pi
\end{eqnarray}
inside the cube implies all the remaining inequalities and then defines a fundamental region $\Delta$. The volume of this region is
\begin{eqnarray}
V=\int_\Delta \prod_{i=1}^r ds_i=\frac 1{\prod_{i=1}^r \tilde n_i} \int_{0\leq y_1+\ldots+y_r\leq \pi, 0\leq y_i \leq \pi}
\prod_{i=1}^r dy_i=\frac {\pi^r}{r!\prod_{i=1}^r \tilde n_i},
\end{eqnarray}
whereas the torus has volume $\pi^r$, so that the number $\nu$ of elementary cells in the cube is:
\begin{eqnarray}
\nu=r!\prod_{i=1}^r \tilde n_i.
\end{eqnarray}
Thus, we can write
\begin{eqnarray}
I=\frac {2^r}{V_F \prod_{i=1}^r \|\alpha_{a_i}\|^2} \frac 1{r!\prod_{i=1}^r \tilde n_i} \int_{Q}
\prod_{a=1}^h |\sin (\vec n_a \cdot \vec s)| \prod_{i=1}^r  ds_i,
\end{eqnarray}
$Q$ being the cube. By setting $2s_i=\zeta_i$ this can also be written as
\begin{eqnarray}\label{dyson}
&& I=\frac {(2\pi)^r}{2^h V_F \prod_{i=1}^r \|\alpha_{a_i}\|^2} \frac 1{r!\prod_{i=1}^r \tilde n_i} J_{\frac 12}, \\
&& J_{\frac 12}= \frac 1{(2\pi)^r}\int_0^{2\pi} d\zeta_1 \ldots \int_0^{2\pi} d\zeta_r \prod_{\alpha \in R} (1- e^{\vec n_\alpha \cdot \vec \zeta})^{\frac 12}.
\end{eqnarray}
Here $J_{\frac 12}$ is a generalized Dyson integral, as conjectured by Macdonald in \cite{mac}, conjecture $2.1"$, for any root system:

\

\

\begin{conj}
For all $s\in \mathbb{C}$ with ${\rm Re}(s)>0$,
\be\label{conjecture}
J_s=\frac 1{(2\pi)^r}\int_0^{2\pi} d\zeta_1 \ldots \int_0^{2\pi} d\zeta_r \prod_{\alpha \in R} (1- e^{\vec n_\alpha \cdot \vec \zeta})^{s} 
=\prod_{i=1}^r \frac {\Gamma (sd_i +1)}{\Gamma(s+1)\Gamma (sd_i -s+1)}.
\ee
\end{conj}

\

\

\noindent
This formula is known as Macdonald's conjecture, in fact it has been proven for all root systems \cite{Opd}. From (\ref{integral}) and
(\ref{dyson}) we get
\begin{eqnarray}
J_{\frac 12} =\frac {2^h V_F r! \prod_{i=1}^r (n_i \|\alpha_{a_i}\|^2)}{\pi^r} \ \frac {{\rm Vol} (\pmb G)}{{\rm Vol} (\pmb H)^2}
\ \frac {|\pmb\Gamma|}{2^r}. \label{int1}
\end{eqnarray}
The last factor is $1$ for $\pmb G$ simply connected.
This formula provides a proof of (\ref{conjecture}) for $s=\frac 12$ and for all the reduced simple lattices.

\section{Arbitrary maximal symmetric embedded subgroups}\label{sec:gms}
As in the previous section we restrict our attention to the case of $\pmb G$ simple.
Here we extend previous results to the general case when $\pmb H$ is not
MCS.
In this case
\begin{eqnarray}
l:= {\rm Rank}(\pmb G/\pmb H) < {\rm Rank} (\pmb G)=r,
\end{eqnarray}
so that the largest possible intersection between the Cartan subalgebra $\pmb {\mathfrak c}$ of $\pmb {\mathfrak g}$ and the complement of
$\pmb {\mathfrak h}$ has dimension $l$. We choose the Cartan subalgebra $\pmb {\mathfrak c}$ just in this way, so that
\begin{eqnarray}
\pmb {\mathfrak c}=\pmb {\mathfrak c}_h \oplus \pmb {\mathfrak c}_p, \qquad \pmb {\mathfrak c}_p:= \pmb {\mathfrak c}\cap \pmb{\mathfrak {p}},
\qquad \dim\pmb {\mathfrak c}_p=l,
\quad \pmb {\mathfrak c}_h\subset \pmb{\mathfrak{h}}.
\end{eqnarray}
Let us fix a basis $k_1,\ldots,k_s$ for $\pmb {\mathfrak c}_h$, $h_1,\ldots, h_l$ for $\pmb {\mathfrak c}_p$, $s+l=r$. Let $\pmb{\mathfrak {k}}$ be the largest Lie
subalgebra of $\pmb{\mathfrak {h}}$ such that $[\pmb{\mathfrak{k}}, \pmb {\mathfrak c}_p]=0$. It is the Lie algebra of the normalizer $\pmb K$ of
$\pmb {\mathfrak c}_p$ in $\pmb H$. Thus, we can write $\pmb {\mathfrak{h}}=:\pmb{\mathfrak {k}}\oplus \pmb{\tilde {\mathfrak {h}}}$ and $\pmb {\mathfrak{p}}=:\pmb{\mathfrak {c}}_p\oplus \pmb{\tilde {\mathfrak {p}}}$, so that:
\begin{eqnarray}
\pmb{\mathfrak{g}}=(\pmb{\mathfrak {k}}\oplus \pmb{\tilde {\mathfrak {h}}}) \oplus (\pmb {\mathfrak c}_p \oplus \pmb{\tilde {\mathfrak {p}}}).
\end{eqnarray}
Since $\pmb{\mathfrak {h}}$ is maximal, we have
\begin{eqnarray}
[\pmb{\mathfrak{h}}, \pmb{\mathfrak{h}}]\subseteq \pmb{\mathfrak{h}}, \qquad
[\pmb{\mathfrak{p}}, \pmb{\mathfrak{p}}]\subseteq \pmb{\mathfrak{h}}, \qquad
[\pmb{\mathfrak{h}}, \pmb{\mathfrak{p}}]\subseteq \pmb{\mathfrak{p}}.
\end{eqnarray}
Moreover, $[\pmb{\mathfrak{k}}, \pmb {\mathfrak c}_p]=0$ implies
\begin{eqnarray}\label{struttura}
[\pmb{\tilde {\mathfrak{h}}}, \pmb {\mathfrak c}_p]\subseteq \pmb{\tilde {\mathfrak {p}}}, \qquad
[\pmb{\tilde {\mathfrak{p}}}, \pmb {\mathfrak c}_p]\subseteq \pmb{\tilde {\mathfrak {h}}}.
\end{eqnarray}
Notice that the roots of $\pmb{\mathfrak g}$ can be divided as follows. Since $\pmb {\mathfrak c}_h$ is the Cartan subalgebra of both
$\pmb{\mathfrak k}$ and $\pmb{\mathfrak h}$,
${\rm Rank}(\pmb{\mathfrak {k}})={\rm Rank}(\pmb{\mathfrak {h}})=s$. We represent the roots as the simultaneous eigenvalues of the operators
$({\rm ad}_{k_1},\ldots, {\rm ad}_{k_s}; {\rm ad}_{h_1},\ldots, {\rm ad}_{h_l})$. The eigenvectors of the roots $\alpha_{\pmb{\mathfrak h}, a}$, $a=1,\ldots, k-s$
($k:=\dim \pmb K$), of $\pmb{\mathfrak {k}}$ are in the complexification
of $\pmb{\mathfrak {k}}$ and thus in the kernel of ${\rm ad}_{h_i}$, $i=1,\ldots,l$: the last $l$ components are zero. Indeed, these are all the
nonvanishing roots with this property, the remaining ones have necessarily non vanishing elements out of the first $s$ ones. We will
call the corresponding roots $\alpha_{\pmb{\mathfrak p}, b}$, $b=1,\ldots, 2q$, where $q$ is the number of positive roots. This root system is not reduced so that each root $\alpha_{\pmb{\mathfrak p}, b}$ has multiplicity $m_b$, and $\sum_{b=1}^q m_b=h-k$. Indeed, these correspond to the non vanishing roots of the
$h_i$. As usual, we can divide all roots in positive and negative, $R=R^+ \oplus R^-$. This
will determine a corresponding decomposition of the restricted root system:
$R_{\pmb{\mathfrak p}}=R_{\pmb{\mathfrak p}}^+ \oplus R_{\pmb{\mathfrak p}}^-$. The main difference w.r.t. the case of a MCS subgroup is that now $R_{\pmb{\mathfrak p}}$ is not a reduced root lattice system and generically each root
$\alpha$ is characterized by a multiplicity $m_\alpha\geq 1$. All such systems are classified in \cite{araki}, see also \cite{Hel}.\\
From now on, we can proceed exactly as in the previous section, by choosing an orthonormal basis of $\pmb{\mathfrak g}$
$B=B_{\pmb K} \cup \{t_1,\ldots,t_{h-k}\} \cup \{h_1,\ldots, h_l\} \cup \{p_1,\ldots, p_{h-k}\},$ where
$B_{\pmb K}=\{k_1,\ldots, k_s, g_1,\ldots, g_{k-s}\}$ is an orthonormal basis for $\pmb{\mathfrak {k}}$, the $t_a$ generate
$\pmb{\tilde {\mathfrak h}}$,
and the $p_b$ generate $\pmb{\tilde {\mathfrak p}}$.
The Euler parameterization for $\pmb G$ is then
\begin{eqnarray}
\pmb G[\vec x; \vec y; \vec z]=e^{\sum_{a=1}^{h-k} x^a t_a} e^{\sum_{i=1}^l y^i h_i} \pmb H[z_1, \ldots ,z_h], \label{euler}
\end{eqnarray}
where $\pmb H$ can be parameterized itself by means of the Euler parameterization, but it is not important here. The range of the $z$ coordinates
must be chosen in such the way to cover the whole subgroup $\pmb H$.
The invariant measure can be computed exactly as in the previous section, giving
\begin{eqnarray}
d\mu_{\pmb G} [\vec x; \vec y; \vec z]=d\mu_{\pmb H} [\vec z] \ d\mu_{\pmb H/\pmb K} [\vec x] \ \prod_{a=1}^{q}
\sin^{m_a} (\vec \alpha_{\pmb{\mathfrak {p}}, a} \cdot \vec y)
\prod_{i=1}^l dy_i,
\end{eqnarray}
where $\vec \alpha_{\pmb{\mathfrak p},a}:= (\alpha_{\pmb{\mathfrak p},a}^1,\ldots, \alpha_{\pmb{\mathfrak p},a}^l)$, $a=1,\ldots, q$ are
the last $l$ components of the positive $\alpha_{\pmb{\mathfrak p}a}$,
corresponding to the eigenvalues of the ${\rm ad}_{h_i}$ only. As before, we can choose a basis of $l$ simple roots
$\vec \alpha_1, \ldots, \vec \alpha_l$ in $R_{\pmb{\mathfrak p}}^+$, to prove
that the range for the coordinates $\vec y$ is given by
\begin{eqnarray}
0\leq \vec \alpha_i \cdot \vec y \leq \pi, \qquad 0\leq \sum_{i=1}^l n_i \vec \alpha_i \cdot \vec y \leq \pi, \label{range}
\end{eqnarray}
where $\sum_{i=1}^l n_i \vec \alpha_i$ is the highest root of the quotient manifold.

\subsection{Further connections with the generalized Dyson integrals}\label{sec:dyson1}
In \cite{mac} Macdonald proposed a generalization of the Dyson integrals extended to not necessarily reduced root lattices. This general conjecture has
been proved by Opdam \cite{Opd} in the form:
\begin{align}
J_{\{k_\alpha\}}^{\pmb{\mathfrak p}}&:=
\frac 1{(2\pi)^r}\int_0^{2\pi} d\zeta_1 \ldots \int_0^{2\pi} d\zeta_r \prod_{\alpha \in R_{\pmb{\mathfrak p}}}
(1- e^{\vec n_\alpha \cdot \vec \zeta})^{k_\alpha} \cr
&=
\prod_{\alpha\in R_{\pmb{\mathfrak p}}^+}\frac{\Gamma(\langle\rho(k),\check{\alpha}\rangle+k_\alpha+\frac 12 k_{\frac \alpha 2}+1)\Gamma(\langle\rho(k),\check{\alpha}\rangle
-k_\alpha-\frac 12 k_{\frac \alpha 2}+1)}{\Gamma(\langle\rho(k),\check{\alpha}\rangle+\frac 12 k_{\frac \alpha 2}+1)\Gamma(\langle\rho(k),\check{\alpha}\rangle
-\frac 12 k_{\frac \alpha 2}+1)},\label{formula}
\end{align}
where $R_{\pmb{\mathfrak p}}$ is a root system, $R_{\pmb{\mathfrak p}}^+$ is a choice of corresponding positive roots,
\be
\rho(k)=\frac 12 \sum_{R_{\pmb{\mathfrak p}}^+} k_\alpha \alpha,
\ee
and $k$ is a Weyl invariant function over $R_{\pmb{\mathfrak p}}$ whose values $k_\alpha$ have positive real part.
For example, the multiplicities $m_\alpha$ select such a function. Finally, $\langle\rho(k),\check{\alpha}\rangle$ indicates the invariant product with the coroot
$\check{\alpha}$.

Repeating the same procedure as in Section \ref{sec:dyson} we get the following formula:
\begin{eqnarray}
J_{\left\{\frac {m_\alpha}2\right\}}^{\pmb{\mathfrak p}}
=\frac {2^{h-k} |\vec \alpha_1 \wedge \ldots \wedge \vec \alpha_l | l! \prod_{i=1}^l n_i}{\pi^l}
\ \frac {{\rm Vol} (\pmb G) {\rm Vol} (\pmb K)}{{\rm Vol} (\pmb H)^2}. \label{int}
\end{eqnarray}
Compared with Theorem 4.1 in \cite{Opd}, with the invariant functions $k_\alpha=m_\alpha/2$, this expression indeed provides the right value for the generalized
Dyson integrals $J_{\{\frac {m_\alpha}2\}}^{\pmb{\mathfrak p}}$, thus a proof of Macdonald's conjecture (\cite{mac}, conjecture 2.3) for $k_\alpha=\frac{m_\alpha}{2}$ and for the lattices associated to all the irreducible symmetric spaces. The ingredients necessary to compute (\ref{int}) are given in Table \ref{tab:fourierFg2}. One then easily
checks, case by case, that formula (\ref{int}) provides the same result as (\ref{formula}).

\section{Euler parameterizations of the simply connected simple Lie groups}\noindent
As an application of our results, we summarize how to realize the generalized Euler parameterization of any simple, simply connected, compact Lie group $\pmb G$ w.r.t. a maximal
symmetrically embedded Lie subgroup $\pmb H$. This is given by expression (\ref{euler}) which we repeat here for convenience:
\begin{eqnarray}
\pmb G[\vec x; \vec y; \vec z]=e^{\sum_{a=1}^{h-k} x^a t_a} e^{\sum_{i=1}^l y^i h_i} \pmb H[z_1, \ldots ,z_h].
\end{eqnarray}
The parameterization of $\pmb H[z_1, \ldots ,z_h]$ can be done inductively in the same way. As we have seen, this is obtain by putting the subgroup $\pmb K$
in evidence so that $z_1, \ldots, z_h$ are chosen in such the way to cover the whole $\pmb H$, the coordinates $x_1, \ldots, x_{h-k}$ have the same
range as $z_1, \ldots, z_{h-k}$. Finally, the range for $y_1, \ldots, y_l$ is specified by (\ref{range}). All possible Euler parameterizations of the simple, simply connected, compact Lie groups
are listed in Table \ref{tab:fourierFg2}.

From the same table one can easily verify that formula (\ref{int}) indeed agrees with (\ref{formula}), thus providing an alternative
proof of the Macdonald-Dyson's conjecture for all simple groups, for the case $k_\alpha=\frac {m_\alpha}2$. The volumes of the groups can be computed as in
\cite{mac1}.

We point out the fact that not all subgroups $\pmb H$ and $\pmb K$ are semisimple but can contain $U(1)$ factors which must be discussed separately. The measure is normalized so that the volume of a $U(1)$ factor is just the length of its period. It is
interesting to notice that such periods can be related to the length of the roots. We will provide a proof of this
fact together with a detailed construction of Table \ref{tab:fourierFg2}, which requires much more space, in a separated publication. Here we limit ourselves to specify the length of the
period for the $U(1)$ factors appearing in the table, after normalizing the long roots of $\pmb G$ to $\sqrt 2$.
They are the following:
\begin{itemize}
\item in the $\mathrm{AIII}_a$ case there is a phase factor in $\pmb H$ with period $T_H=2\pi\sqrt{\frac{p+q}{pq}}$, whereas $\pmb K$ contains $p$ phase factors with periods
$T_i=\frac{2\pi}{i}\sqrt{2i(i+1)}$, for $i=1,\cdots,p-1$ and $T_p=2\pi\sqrt{\frac{2p(p+q)}{q-p}}$;
\item in the $\mathrm{AIII}_b$ case there is a phase factor in $\pmb H$ with period $T_{H}=2\pi\sqrt{\frac 2p}$, and $p-1$ phase factors in $\pmb K$ with periods $T_i=\frac{2\pi}{i}\sqrt{2i(i+1)}$, for $i=1,\cdots,p-1$;
\item in the $\mathrm{AIV}$ case there is a phase factor in $\pmb H$ with period $T_{H}=2\pi\sqrt{\frac{n+1}{n}}$, and a phase factor in $\pmb K$ with period $T_K=2\pi\sqrt{\frac{n+1}{2(n-1)}}$;
\item in the $\mathrm{CI}$ case the phase factor in $\pmb H$ has period $T_H=2\pi\sqrt{2n}$;
\item in the $\mathrm{DI}_b$ case the phase factor in $\pmb K$ has period $T_K=4\pi$;
\item in the $\mathrm{DIII}_a$ case the subgroup $\pmb H$ is $U(2n+1)\simeq SU(2n+1)\times U(1)/\mathbb{Z}_{2n+1}$ and the period of the
phase factor is $T_H=4\pi\sqrt{2n+1}$, whereas the phase factor in $\pmb K$ has period $T_K=4\pi$;
\item in the $\mathrm{DIII}_b$ case the subgroup $\pmb H$ is $U(2n)\simeq SU(2n)\times U(1)/\mathbb{Z}_{2n}$ and the period of the phase factor
is $T_H=4\pi\sqrt{2n}$;
\item in the $\mathrm{EII}$ case the two phase factors in $\pmb K$ have periods $T_{K_1}=4\pi$ and $T_{K_2}=4\pi\sqrt{3}$;
\item in the $\mathrm{EIII}$ case the periods of the phase factors in $\pmb H$ and in $\pmb K$ are $T_H=4\pi\sqrt{3}$ and $T_K=4\pi\sqrt{3}$;
\item in the $\mathrm{EVII}$ case the period of the phase factor in $\pmb H$ is $T_H=2\pi\sqrt{\frac{3}{2}}$.
\end{itemize}
Moreover, there are some particular cases that must be considered separately in the table, so that we list them apart:
\begin{itemize}
\item $\mathrm{AI}$: for $n=1$, $\pmb H=SO(2)$ with period $T=4\pi$ and obviously $\alpha_h$ is not defined; \\
\phantom{$AI$: }for $n=2$, $\pmb H=SO(3)$ which has only the short root, so that $|\alpha_G|/|\alpha_H|=2$;
\item $\mathrm{BI}_a$: for $n=2$, $\pmb H=SO(2)\times SO(3)$ whose phase factor has period $T=4\pi$ and $|\alpha_G|/|\alpha_H|=\sqrt{2}$; \\
\phantom{$BI_a$: }for $n=3$, $\pmb H=SO(3)\times SO(4)$ and the ratios of the root lengths are $|\alpha_G|/|\alpha_{SO(3)}|=\sqrt{2}$ and $|\alpha_G|/|\alpha_{SO(4)}|=1$;
\item $\mathrm{BI}_b$: for $p=2$ and $q>3$, $\pmb H=SO(2)\times SO(q)$ whose phase factor has period $T=4\pi$ and the ratio of the root lengths is $|\alpha_G|/|\alpha_{SO(q)}=1$; \\
\phantom{$BI_b$: }for $p=3$ and $q>5$, $\pmb H=SO(3)\times SO(q)$ and the ratios of the root lengths are $|\alpha_G|/|\alpha_{SO(3)}|=\sqrt{2}$ and $|\alpha_G|/|\alpha_{SO(q)}|=1$;
\item $\mathrm{BII}$: for $n=1$ is the same as $AI$ for $n=1$;
\item $\mathrm{DI}_a$: for $n=2$, $\pmb H=SO(2)\times SO(2)$ whose phase factors have both period $T=4\pi$;
\item $\mathrm{DI}_b$: for $n=3$, $\pmb H=SO(2)\times SO(4)$ whose phase factor has period $T=4\pi$;
\item $\mathrm{DI}_c$: for p=2 $q>3$, $\pmb H=SO(2)\times SO(q)$ whose phase factor has period $T=4\pi$; \\
\phantom{$DI_c$: }for $p=3$ and $q>4$, $\pmb H=SO(3)\times SO(q)$ and the ratios of the root lengths are $|\alpha_G|/|\alpha_{SO(3)}|=\sqrt{2}$ and $|\alpha_G|/|\alpha_{SO(q)}|=1$;
\item $\mathrm{DII}$: for $n=2$, $\pmb H=SO(3)$ and the ratio of the root lengths is $|\alpha_G|/|\alpha_{SO(3)}|=\sqrt{2}$.
\end{itemize}
In this list the unspecified data can be read from Table \ref{tab:fourierFg2}.\\

\begin{landscape}
\begin{table}[hbtp]
\begin{center}
\resizebox*{1.3\textheight}{!}{
\begin{tabular}{ccccccccccccccc}
\toprule \pmb{Label} & $\rm{\pmb{G}}_{c}$ & $\rm{\pmb{G}}_{nc}$ & dim(\pmb G) & $\pmb Z$ & MCS \pmb H & dim(\pmb H) & $\pmb{\Lambda_{G/H}}$ &
$\pmb{(n_1,\ldots, n_r)}$  & $\pmb{|\alpha_G|/|\alpha_H|}$ & $\pmb{\vec{m}_\lambda,\, \vec{m}_{2\lambda}}$ & $\pmb{|\alpha_G|/|\alpha_{G/H}|}$ & $\pmb K$ & $\pmb{|\alpha_H|/|\alpha_K|}$ & ${\pmb \rho}$\\
\midrule
AI &\pmb{SU}(n+1) &$\pmb{SL}(n+1,\mathbb{R})$ & $n^2+2n$ & $\pmb{\mathbb Z}_{n+1}$ & \pmb{SO}(n+1)/$\pmb{\mathbb Z}_{2}$ & $n(n+1)/2$ & $A_n \,(n\geq 1)$ &(1,1,\ldots,1)  & $\sqrt{2}$
& (1), (0) & 1 & $\pmb{\mathbb Z}_2^n$ & - & $V_{\lambda_1}$ \\[0.2cm]
AII & $\pmb{\rm{SU}}(2n)$ & $\pmb{\rm{SU}}^*(2n)$ & $4n^2-1 $  & $\pmb{\mathbb Z}_{2n}$ & \pmb{USp}(2n) & $2n^2+n$ & $A_{n-1}\,(n>1)$ & (1,1,\ldots,1)  & 1
& (4), (0) & $\sqrt{2}$ & $\pmb{\rm{SU}}(2)^{n}$ & 1 & $V_{\lambda_1}$ \\[0.2cm]
$\rm{AIII}_a$ & \pmb{SU}(p+q) & \pmb{SU}(p,q)& $(p+q)^2-1$ &  $\pmb{\mathbb Z}_{p+q}$ & $\pmb{S}(\pmb{U}(p)\times \pmb{U}(q))$ & $p^2+q^2-1$ & $B_p\, (1<p<q)$ &
(1,2,\ldots,2) & 1 & $2(1,q-p), (0,1)$ & 1 & $S(\pmb{\rm{U}}(1)^p\times\pmb{\rm{U}}(q-p))$ & 1 & $V_{\lambda_1}$ \\[0.2cm]
$\rm{AIII}_b$ &\pmb{SU}(2p)& \pmb{SU}(p,p) & $4p^2-1$ & $\pmb{\mathbb Z}_{2p}$ & $\pmb{S}(\pmb{U}(p)\times \pmb{U}(p))$ & $2p^2-1$ & $C_p\,(p>1)$ &
$(2,2,\ldots,2,1)$ & 1 & (1,2), (0,0) & 1 & $\pmb{S}(\pmb{U}(1)^p)=\pmb{SO}(2)^{p-1}\times \pmb{\mathbb Z}_2$ & 1 & $V_{\lambda_1}$ \\[0.2cm]
AIV & $\pmb{\rm{SU}}(n+1)$ & $\pmb{\rm{SU}}(1,n)$ &$n^2+2n$ & $\pmb{\mathbb Z}_{n+1}$ & $\pmb{S}(\pmb{U}(1)\times \pmb{U}(n))$  & $n^2$ & $A_1$ & (2) & 1 &
(2n-2), (1) & 1 & $\pmb{S}(\pmb{U}(n-1)\times \pmb{U}(1)) $ & 1 & $V_{\lambda_1}$ \\[0.2cm]
$\rm{BI_a}$ & $\pmb{\rm{SO}}(2n+1)$ &$\pmb{\rm{SO}}_0(n,n+1)$ & $2n^2+n$  & $\pmb {\mathbb Z}_2$ & \pmb{SO}(n)$\times$\pmb{SO}(n+1) & $n^2$ & $B_n\,(n> 3)$
&(1,2,\ldots,2) & 1 & (1,1),\, (0,0) & 1 & $\pmb{\mathbb Z}_2^n$ & - & $V_{\lambda_n}$ \\[0.2cm]
$\rm{BI_b}$ & $\pmb{\rm{SO}}(p+q)=\pmb{\rm{SO}}(2n+1)$ &$\pmb{\rm{SO}}_0(p,q)$ & (p+q)(p+q-1)/2 & $\pmb {\mathbb Z}_2$ & \pmb{SO}(p)$\times$\pmb{SO}(q) &
$p(p-1)/2+q(q-1)/2$ & $B_p\, (1<p<n)$ & (1,2,\ldots,2)  & 1 & (1,2(n-p)+1), (0,0) &1 & $\pmb{SO}(q-p)\times \pmb {\mathbb Z}_2 \ltimes \pmb {\mathbb Z}_2^{p} $
& 1 & $V_{\lambda_n}$ \\[0.2cm]
BII & $\pmb{\rm{SO}}(2n+1)$ &$\pmb{\rm{SO}}_0(1,2n)$ & $2n^2+n$ & $\pmb {\mathbb Z}_2$ & \pmb{SO}(2n) & n(2n-1) & $A_1$ & (1) & 1 & (2n-1), (0) & $\sqrt{2}$ &
$\pmb{SO}(2n-1)$ & 1 & $V_{\lambda_n}$ \\[0.2cm]
CI & \pmb{USp}(2n) & $\pmb{\rm{Sp}}(2n,\mathbb{R})$ & $2n^2+n$ & $\pmb {\mathbb Z}_2$ & \pmb{U}(n) & $n^2$ & $C_{n}\,(n\geq 3)$ &(2,2,\ldots,2,1) & $\sqrt{2}$ &
(1,1),\, (0,0) & 1 & $\pmb {\mathbb Z}_2^n$ & - & $V_{\lambda_1}$ \\[0.2cm]
$CII_a$ & \pmb{USp}(2p+2q)=\pmb{USp}(2n) & \pmb{USp}(2p,2q) &  $(p+q)(2p+2q+1)$ & $\pmb {\mathbb Z}_2$ & \pmb{USp}(2p)$\times$\pmb{USp}(2q) & $2p^2+p+2q^2+q$ &
$B_p\, (1\leq p\leq(n-1)/2)$ & (2,2,\ldots,2) & 1 & (4,4n-8p), (0,3) &1 & $\pmb{USp}(2q-2p)\times \pmb{SU}(2)^p$ & 1 & $V_{\lambda_1}$ \\[0.2cm]
$CII_b$ & \pmb{USp}(4n) & \pmb{USp}(2n,2n) & $8n^2+2n$ & $\pmb {\mathbb Z}_2$ & \pmb{USp}(2n)$\times$\pmb{USp}(2n) &$4n^2+2n$ & $C_n$ & $(2,2,\ldots,2,1)$ & 1
& (3,4), (0,0) & $\sqrt{2}$ & $\pmb{SU}(2)^n$ & $\sqrt{2}$ & $V_{\lambda_1}$\\[0.2cm]
$DI_a$ &\pmb{SO}(2n) & \pmb{SO}(n,n)  & $n(2n-1)$ & $Z$ &\pmb{SO}(n)$\times$\pmb{SO}(n) & $n(n-1)$ & $D_{n}\,(n > 1)$ &(1,2,\ldots,2,1,1) & 1 & (1), (0) &
1 & $\pmb {\mathbb Z}_2^{n}$ & - & $V_{\lambda_n}$ \\[0.2cm]
$DI_b$ &\pmb{SO}(2n)  & \pmb{SO}(n-1,n+1) & $n(2n-1)$ & $Z$  &\pmb{SO}(n-1)$\times$\pmb{SO}(n+1) & $n^2-n+1$ & $B_{n-1} \,(n> 2)$ &(1,2,\ldots,2) & 1 &
(1,2), (0,0) & 1 &$\pmb{U}(1)\times \pmb {\mathbb Z}_2 \ltimes \pmb {\mathbb Z}_2^{n-1} $ & 1 & $V_{\lambda_n}$\\[0.2cm]
$DI_c$ &\pmb{SO}(p+q)=\pmb{SO}(2n)  & \pmb{SO}(p,q) & (p+q)(p+q-1)/2 & $Z$  &\pmb{SO}(p)$\times$\pmb{SO}(q) & $p(p-1)/2+q(q-1)/2$ & $B_{p} \,(1<p<n-1)$
&(1,2,\ldots,2) & 1 & (1,2(n-p)), (0,0) & 1 & $\pmb{SO}(q-p)\times \pmb {\mathbb Z}_2 \ltimes \pmb {\mathbb Z}_2^p $ & 1 & $V_{\lambda_n}$ \\[0.2cm]
$DII$ &\pmb{SO}(2n) & \pmb{SO}(1,2n-1)  & $n(2n-1)$ & $Z$  &\pmb{SO}(2n-1) & $(2n-1)(n-1)$ & $A_1$ &(1) & 1 & (2n-2), (0) & 1 & \pmb{SO}(2n-2) & 1 & $V_{\lambda_n}$ \\[0.2cm]
$DIII_a$ &\pmb{SO}(4n+2)  & $\pmb{SO}^*(4n+2)$ & $(2n+1)(4n+1)$ & $\pmb {\mathbb Z}_4$ &\pmb{U}(2n+1) & $(2n+1)^2$ & $B_n \,(n\geq 2)$ & (2,2,\ldots,2) & 1 &
(4,4), (0,1) &1&$\pmb{SU}(2)^n\times\pmb{SO}(2)$&2 & $V_{\lambda_{2n+1}}$ \\[0.2cm]
$DIII_b$ &\pmb{SO}(4n) & $\pmb{SO}^*(4n)$  & $2n(4n-1)$ & $\pmb {\mathbb Z}_2\times \pmb {\mathbb Z}_2$ &\pmb{U}(2n) & $4n^2$ & $C_n\,(n \geq 2)$
&(2,2,\ldots,2,1) & 1 & (1,4), (0,0) & 1 & $\pmb{SU}(2)^n$ & 1 & $V_{\lambda_{2n}}$ \\[0.2cm]
EI &\pmb{E}$_{6(-78)}$ & \pmb{E}$_{6(6)}$ & 78 & $\pmb{\mathbb Z}_3$ & \pmb{USp}(8)/$\pmb{\mathbb{Z}_2}$ & 36 & $E_6$ & (1,2,2,3,2,1) & 1 &(1), (0) & 1 &
$\pmb{\mathbb Z}_2^6$ & - & $V_{\lambda_1}$ \\[0.2cm]
EII &\pmb{E}$_{6(-78)}$ & \pmb{E}$_{6(2)}$ & 78 & $\pmb{\mathbb Z}_3$ &  (\pmb{USp}(2)$\times$\pmb{SU}(6))/$\pmb{\mathbb{Z}}_2$ & 38 & $F_4$ & (2,3,4,2) & 1
& (1,2), (0,0) &1 & $\pmb{U}(1)^2\times \pmb{\mathbb{Z}}_2$ & 1 & $V_{\lambda_1}$ \\[0.2cm]
EIII &\pmb{E}$_{6(-78)}$ & \pmb{E}$_{6(-14)}$ & 78 & $\pmb{\mathbb Z}_3$ & (\pmb{U}(1)$\times$\pmb{SO}(10))/$\pmb{\mathbb{Z}}_4$ & 46 & $B_2$ & (2,2) & 1 &
(6,8), (0,1) & 1 & (\pmb{SO}(6)$\times$\pmb{U}(1))/$\pmb{\mathbb{Z}}_2$ & 1 & $V_{\lambda_1}$ \\[0.2cm]
EIV &\pmb{E}$_{6(-78)}$ & \pmb{E}$_{6(-26)}$ & 78 & $\pmb{\mathbb Z}_3$ & $\pmb{F}_4$ & 52 & $A_2$ & (1,1) & 1 & (8), (0) & $\sqrt{2}$ & $\pmb{SO}(8)$ & 1 & $V_{\lambda_1}$ \\[0.2cm]
EV &\pmb{E}$_{7(-133)}$ & \pmb{E}$_{7(7)}$ & 133 & $\pmb{\mathbb Z}_2$ & \pmb{SU}(8)/$\pmb{\mathbb{Z}}_2$ & 63 & $E_7$ & (2,2,3,4,3,2,1) & 1 & (1), (0) & 1 &
$\pmb{\mathbb Z}_2^7$ & - & $V_{\lambda_6}$\\[0.2cm]
EVI &\pmb{E}$_{7(-133)}$ & \pmb{E}$_{7(-5)}$ & 133 & $\pmb{\mathbb Z}_2$ & $(\pmb{SU}(2)\times \pmb{SO}(12))/\pmb{\mathbb{Z}}_2$  & 69 & $F_4$ & (2,3,4,2) & 1
& (1,4), (0,0) & 1 & $\pmb{SU}(2)^3\times\pmb{\mathbb{Z}}_2\times\pmb{\mathbb{Z}}_2$ & 1 & $V_{\lambda_6}$ \\[0.2cm]
EVII &\pmb{E}$_{7(-133)}$ & \pmb{E}$_{7(-25)}$ & 133 & $\pmb{\mathbb Z}_2$ & $(\pmb{U}(1)\times E_6)/\pmb{\mathbb{Z}}_3$ & 79 & $C_3$ & (2,2,1) & 1 &(1,8),
(0,0) & $\sqrt{2}$ & $\pmb{SO}(8)$ & 1& $V_{\lambda_6}$ \\[0.2cm]
EVIII &\pmb{E}$_{8(-248)}$ & \pmb{E}$_{8(8)}$ & 248 & \pmb 1 & \pmb{Ss}(16) & 120 & $E_8$ &(2,3,4,6,5,4,3,2) &  1 & (1), (0) & 1 &
$\pmb{\mathbb Z}_2^8$ & - & Any \\[0.2cm]
EIX &\pmb{E}$_{8(-248)}$ & \pmb{E}$_{8(-24)}$ & 248 & \pmb 1 & $(\pmb{SU}(2)\times E_7)/\pmb{\mathbb {Z}}_2$ & 136 & $F_4$ & (2,3,4,2) & 1 & (1,8), (0,0) & 1 &
$\pmb{SO}(8)\times \pmb{\mathbb{Z}}_2\times \pmb{\mathbb{Z}}_2$ &1 & Any \\[0.2cm]
FI &\pmb{F}$_{4(-52)}$ & \pmb{F}$_{4(4)}$ & 52 & \pmb 1 & (\pmb{USp}(6)$\times$\pmb{USp}(2))/$\pmb{\mathbb{Z}}_2$ & 24 & $F_4$ & (2,3,4,2) & 1 & (1,1), (0,0) & 1 &
$\pmb{\mathbb Z}_2^4$ & - & Any \\[0.2cm]
FII &\pmb{F}$_{4(-52)}$ & \pmb{F}$_{4(-20)}$ & 52 & \pmb 1 & \pmb{SO}(9) & 36 & $A_1$ &(2) & 1 & (8), (7) & $2\sqrt{2}$ & \pmb{SO}(7) & 1 & Any \\[0.2cm]
G &\pmb{G}$_{2(-14)}$ & \pmb{G}$_{2(2)}$ & 14 & \pmb 1 &\pmb{SO}(4)$/\pmb{\mathbb{Z}}_2$ & 6 & $G_2$ &(3,2) & 1 & (1,1), (0,0) & 1 & $\pmb{\mathbb Z}_2^2$ & - & Any \\
\bottomrule
\end{tabular}
}
\caption{Maximal symmetrically embedded proper subgroups for the compact simple Lie algebras. Note that we are referring
to the universal coverings, so that $\pmb{SO}(3)\simeq \pmb{SU}(2)$, $\pmb{USp}(4)\simeq \pmb{SO}(5)$, $\pmb{USp}(2)\simeq \pmb{SU}(2)$,
$\pmb{SO}(6)\simeq \pmb{SU}(4)$, $\pmb{SO}(4)\simeq \pmb{SU}(2)\times \pmb{SU}(2)$, and $\pmb {SO}(n)\simeq \pmb {Spin}(n)$; in ${\rm EVII}$ $\pmb{Ss}(16)\simeq \pmb{SO}(16)/\pmb{\mathbb Z}_2$ is a semispin group.
In the second column we indicate the compact form associated to the real form listed in the third column. $Z$ indicates the center of the compact
form. In particular, $Z$ is ${\mathbb Z}_4$ if the dimension of the spin group is $4k+2$ and ${\mathbb Z}_2\times {\mathbb Z}_2$ if the dimension is $4k$.
In the column $\Lambda_{G/H}$ we indicate the reduced form of the root system associated to the symmetric space $G/H$. However, these in general can contain
also double roots. Notice that the rank of the reduced system gives the rank of the symmetric space. The quotients $|\alpha_G|/|\alpha_H|$,
$|\alpha_G|/|\alpha_{G/H}|$ and $|\alpha_H|/|\alpha_K|$ indicates the ratio of the long roots (including eventual double roots) of the indicated
root systems. $(n_1,\ldots, n_r)$ are the coefficients of the highest root of the root system for the symmetric manifold.
$m_\lambda=(m_{\lambda_l},m_{\lambda_s})$ and $m_{2\lambda}=(m_{2\lambda_l},m_{2\lambda_s})$ indicate the multiplicities of the roots of the
reduced lattice and of the double roots respectively, where $l$ and $s$ denote long and short respectively. In the last column $\rho$ denotes a choice for a $G$-faithful representation of the algebra; $V_{\lambda_i}$ means the fundamental representation associated to the i-th weight in the corresponding Dynkin diagram.
}
\label{tab:fourierFg2}
\end{center}
\end{table}
\end{landscape}
\noindent

\section*{Acknowledgments.}
We would like to acknowledge Alessio Marrani and Bert Van Geemen for useful comments. We thank Stefano Pigola 
for helpful explanations. We also thank Daniel Bump for relevant suggestions and remarks.


\end{document}